\documentclass[12pt]{article}
\usepackage{amssymb}
\usepackage{graphicx}
\usepackage{latexsym}

\def\part#1{\frac{\partial\phantom{q}}{\partial#1}}

\newenvironment{rmk}{\begin{trivlist}\item[]{\bf Remark:} }
{\end{trivlist}}

\newtheorem{thm}{Theorem}

\newtheorem{prp}[thm]{Proposition}

\def\End{\mathop{\rm End}\nolimits}

\def\ker{\mathop{\rm ker}\nolimits}
\def\coker{\mathop{\rm coker}\nolimits}

\def\Pic{\mathop{\rm Pic}\nolimits}

\def\rk{\mathop{\rm rk}\nolimits}
\def\tr{\mathop{\rm tr}\nolimits}

\def\ad{\mathop{\rm ad}\nolimits}

\def\Jac{\mathop{\rm Jac}\nolimits}

\def\Imag{\mathop{\rm Im}\nolimits}

\newcommand{\R}{\mathbf{R}}
\newcommand{\C}{\mathbf{C}}
\newcommand{\K}{\mathbf{H}}
\newcommand{\Z}{\mathbf{Z}}

\begin{document}
\title{Integrable systems and Special K\"ahler metrics}
\author{Nigel Hitchin\\{Mathematical Institute,
Woodstock Road,
Oxford, OX2 6GG}\\{hitchin@maths.ox.ac.uk}}

\maketitle

\section{Introduction}
The base of an algebraic completely integrable system has a natural differential geometric structure called a Special K\"ahler structure, describing in a rather different way the variation of the period matrix of the abelian varieties which form the fibres. The moduli space ${\mathcal M}$ of Higgs bundles on a curve $\Sigma$ of genus $g>1$ gives a particularly interesting example -- the so-called Hitchin system. In \cite{BH} the authors give a simple formula for the K\"ahler potential of the Special K\"ahler structure in this case. Here  we give first a short survey of this situation and afterwards extend it to cover some distinguished integrable subsystems.

A Higgs bundle is a pair $(V,\Phi)$ consisting of a holomorphic vector bundle $V$ of rank $n$ and a holomorphic section (the Higgs field) of $\End V\otimes K$, where $K$ is the canonical bundle,  satisfying a joint stability condition. The integrable system arises from taking the curve $S$ defined by  $\det(x-\Phi) = x^n+a_1x^{n-1}+\cdots +a_n=0$, and realizing  $(V,\Phi)$ by the direct image construction from a line bundle on $S$. The spectral curve sits naturally in the total space of the cotangent bundle $T^*\Sigma$ and  the formula in \cite{BH} is
$$K=\Imag \frac{1}{4}\int_S\theta\wedge\bar\theta$$
where $\theta$ is the canonical one-form on the cotangent bundle. 

The assumption here is that $S$ is smooth. However, the formula makes sense for a singular spectral curve and our aim is to show that indeed it holds on the subintegrable systems which come from fixing the singularity type -- here we restrict to ordinary singularities of $S$, points of multiplicity $m$ with $m$ distinct tangents.  For a smooth curve $S$ the fibre of the integrable system is its Jacobian.  For a singular curve it is  the Jacobian of the normalization. 

We end with a discussion of examples and links to the asymptotics of the natural  hyperk\"ahler metric on the moduli space ${\mathcal M}$.

\section{Integrable systems}

\subsection{The geometry of integrable systems}
By an integrable system we mean a symplectic manifold  $(M,\omega)$ of dimension $2n$ and a proper map $h:M\rightarrow B$ to an $n$-dimensional manifold $B$ with Lagrangian fibres. The connected component of the fibre over a regular value is a torus, and functions on $B$ composed with $h$ Poisson commute. For the moment we work over the real numbers.

A tangent vector $X\in T_aB$ lifts canonically to a section of the normal bundle of the fibre $M_a$ -- the first order deformation of $M_a \subset M$ -- and locally we can extend this to a vector field $\tilde X$ along $M_a$, well-defined up to the addition of a vector field $Y$ tangent to $M_a$. But $M_a$ is Lagrangian, so the interior product $i_Y\omega$ restricts to zero and therefore $i_{\tilde X}\omega$ is a well-defined global 1-form  on  the fibre. It is in fact a {\it closed} 1-form. 
The de Rham cohomology class $[i_{\tilde X}\omega]$ now provides a natural isomorphism $T_aB\cong H^1(M_a,\R)$.

If $U\subset B$ is a contractible neighbourhood of $a$ then homotopy invariance means the restriction map $H^1(h^{-1}(U),\R)\rightarrow H^1(M_a,\R)$ is an isomorphism. If $a,b\in U$ then this gives a natural isomorphism $T_aB\cong T_bB$. Globally this is a flat connection on the tangent bundle of $B$. Since $\omega$ vanishes on $M_a$ its cohomology class in $H^2(h^{-1}(U),\R)$ is zero and there we can write $\omega=d\theta$. Choosing a basis $C_i\in H_1(M_a,\Z)\cong \Z^n$ one obtains local flat coordinates 
$$x_i=\int_{C_i}\theta.$$
It follows that the connection $\nabla$ is flat and torsion-free.

Now suppose that $M$ is a holomorphic symplectic manifold and $h:M\rightarrow B$ is holomorphic. Then the real and imaginary parts of the holomorphic symplectic form $\omega^c$ are real symplectic forms.  Let $\omega$ be the real part. Suppose in addition that $M$ has a K\"ahler form $\Omega$. Then the cohomology class $[\Omega]\in H^2(M,\R)$ defines a real symplectic form on the base:
$$(X,Y)=\int_{M_a}i_{\tilde X}\omega\wedge i_{\tilde Y}\omega\wedge \Omega^{n-1}.$$
This is just the cup product of cohomology classes and so is constant with respect to the flat connection. Hence the base $B$ is equipped with a flat symplectic connection. 

The complex structure on the base is not flat. It is an endomorphism $I:T\rightarrow T$ such that $I^2=-1$ or equivalently a section of $T^*\otimes T$. This is a 1-form with values in the flat vector bundle $T$ and the condition it satisfies is $d_{\nabla}I=0$ where $d_{\nabla}$ is the flat de Rham differential operator $d_{\nabla}:\Omega^1(B,T)\rightarrow \Omega^2(B,T)$. Locally $I=d_{\nabla}X$ for a Hamiltonian vector field.

\begin{rmk}
The above data -- a flat torsion-free symplectic connection and a complex structure $I$, compatible with the symplectic form and satisfying  $d_{\nabla}I=0$  -- is the definition of a {\it Special K\"ahler metric} on $B$ in the approach of D.Freed \cite{DF}, directed towards  mathematicians. The original concept appeared in the physics literature in the context of supersymmetry \cite{WP},\cite{G}.
\end{rmk}

We could have chosen the imaginary part of $\omega^c$ or the real part of $e^{i\theta}\omega^c$ to get other flat connections and this is all part of the rich geometry of Special K\"ahler stuctures but we shall restrict ourselves here to one aspect, which is a convenient description of the metric given in \cite{DF} in terms of a K\"ahler potential.  

In the holomorphic situation the fibre $M_a$ is a complex torus and the K\"ahler class of $\Omega$ means we can choose a symplectic basis $A_i,B_i$ of $H_1(M_a,\R)$. Let $x_i,y_i$, $1\le i\le n$ be the real flat coordinates corresponding to this basis, then they are the real parts of holomorphic coordinates $z_i$ and $w_i$ and a K\"ahler potential is given by 
$$K=\frac{1}{2} \Imag\sum_i w_i\bar z_i.$$
If $\omega^c=d\theta^c$ and $\alpha_i,\beta_i$ denote the dual basis of $H^1(M_x,\R)$ then $[\theta^c]=\sum_iz_i\alpha_i+w_i\beta_i$ and the above formula can be written as 
\begin{equation}
K=\Imag \int_{M_a}\theta^c\wedge \bar\theta^c\wedge \Omega^{n-1}.
\label{pot1}
\end{equation} 
\subsection{Higgs bundles} 
The particular holomorphic integrable system we want to apply this formula to is where $M={\mathcal M}$ is the moduli space of Higgs bundles on a Riemann surface $\Sigma$ of genus $g>1$.  A Higgs bundle is a rank $n$ vector bundle $V$ together with a holomorphic section $\Phi$ of $\End V\otimes K$ where $K$ is the canonical bundle. A stability condition for the pair gives rise to a well-behaved moduli space which is a holomorphic symplectic manifold. 

We can represent a Higgs bundle as a pair $(A,\Phi)$ where $A$ denotes a $\bar\partial$-operator 
$\bar\partial_A:\Omega^0(\Sigma, V)\rightarrow \Omega^{01}(\Sigma, V)$ on the rank $n$ holomorphic vector bundle $V$  
and $\Phi$, the Higgs field,   is a section of $\End V\otimes K$ such that 
 $\bar\partial_A\Phi=0$.

In this Dolbeault form a tangent vector to the moduli space ${\mathcal M}$ is given  by a pair $(\dot A,\dot \Phi)$ 
satisfying  $\bar\partial_A\dot\Phi+[\dot A,\Phi]=0$ modulo the action of a $C^{\infty}$ endomorphism $\psi$ defining $(\dot A,\dot\Phi)=(\bar\partial_A\psi,-[\psi,\Phi])$. The symplectic form is then 
\begin{equation}
\omega^c((\dot A_1,\dot\Phi_1),(\dot A_2,\dot\Phi_2))=\int_{\Sigma}\tr(\dot A_1\dot\Phi_2-\dot A_2\dot\Phi_1).
\label{symp}
\end{equation}

The characteristic polynomial of $\Phi$ gives 
$$\det(x-\Phi)=x^n+a_1x^{n-1}+\cdots + a_n$$
where the coefficient $a_m$ is a holomorphic section of $K^m$. Since the characteristic polynomial is conjugation-invariant these coefficients give a well-defined map 
$$h:{\mathcal M}\rightarrow \bigoplus_{i=1}^n H^0(\Sigma, K^i).$$
This defines the integrable system. Later it will be useful to use the invariant polynomials $\tr(\Phi^m)$ instead of the coefficients $a_i$ to define an equivalent map $h$.

A point $a=(a_1,\dots, a_n)$ in the base ${\mathcal B}$ defines an algebraic curve $\det(x-\Phi)=0$, the spectral curve $S$. This is an $n$-fold covering $\pi:S\rightarrow \Sigma$ and a line bundle $L$ on $S$ defines a rank $n$ bundle $V$ on $\Sigma$ by the direct image $\pi_*L$: over an open set $U\subset \Sigma$, 
$H^0(U, V) = H^0(U,\pi_*L)=H^0(\pi^{-1}(U), L)$ by definition. The indeterminate $x$ in the characteristic polynomial is a section of $\pi^*K$ on $S$ and then $x:H^0(\pi^{-1}(U), L)\rightarrow H^0(\pi^{-1}(U), L\pi^*K)$ defines a Higgs field $\Phi$ on $V$. In this way we can identify the fibre over a point $a$ which defines a smooth spectral curve (a regular value of $h$) with the Picard variety $\Pic^d(S)$ for a fixed $d$ depending on the degree of $V$. This is a torsor for the Jacobian of $S$ and is a complex torus. 

The family of $\bar\partial$-operators $\bar\partial_A$ defines a determinant line bundle whose cohomology class is given by a K\"ahler (actually hyperk\"ahler) form $\Omega$  and restricts to the theta class on the Jacobian \cite{BNR}. So we have all the data for a Special K\"ahler metric on the base ${\mathcal B}$.

We have one more feature which is the $\C^*$-action defined by $(V,\Phi)\mapsto (V,\lambda\Phi)$. From (\ref{symp}) it acts on the symplectic form $\omega^c$ by rescaling and so if $X$ is the holomorphic vector field on ${\mathcal M}$ generated by the infinitesimal action we have ${\mathcal L}_X\omega^c=\omega^c$ or $d(i_X\omega^c)+i_Xd\omega^c=\omega^c$ which, since $\omega^c$ is closed, means that $\omega^c=d(i_X\omega^c).$
Hence from (\ref{pot1}) we have a formula for the K\"ahler potential of the Special K\"ahler metric 
$$K=\Imag \int_{{\mathcal M}_a}i_X\omega^c\wedge \overline {i_X\omega^c}\wedge \Omega^{n-1}.$$

Note the geometric interpretation of the 1-form $i_X\omega^c$ on a fibre ${\mathcal M}_a$. If we regard it as a tangent vector on the base, or a first order deformation of the fibre, it is simply the deformation given by the action of $\C^*$.

\subsection{Spectral curves} 
The base ${\mathcal B}$ of the integrable system now has two interpretations: the parameter space of Lagrangian fibres in ${\mathcal M}$ and the parameter space of spectral curves. Now, as shown in \cite{Hit1}, a special K\"ahler structure exists on the moduli space of deformations of a compact complex Lagrangian in a complex symplectic K\"ahler manifold, and not just the fibres of an integrable system. The procedure is the same -- a section of the normal bundle defines a holomorphic 1-form and its periods define the Special K\"ahler structure. 

A special case consists of the spectral curves themselves. The equation $x^n+a_1x^{n-1}+\cdots + a_n=0$ defines the zero set of a section of $\pi^*K^n$ on the total space $\vert K\vert$ of $K$ where $x$ is the tautological section of $\pi^*K$ on $\vert K\vert$. Thus $x$ embeds $S$ as a curve in the complex surface $\vert K\vert$,  the cotangent bundle $T^*\Sigma$ of $\Sigma$, and the covering is given by the restriction of the projection $\pi:T^*\Sigma \rightarrow \Sigma$.  The  cotangent bundle is canonically symplectic and any curve is Lagrangian. Moreover, when $n=1$  the factor $\Omega^{n-1}$ disappears in the definition of  the symplectic form on the moduli space, so we do not require a K\"ahler form on the surface.

All sections of  $\pi^*K^n$  on $\vert K\vert$ are of the form $a_0x^n+a_1x^{n-1}+\cdots + a_n$ and compact zero sets require $a_0\ne 0$, thus any deformation is described by the equation of a spectral curve. In particular since  the normal bundle to a curve in a symplectic surface is the canonical bundle $K_S\cong \pi^*K^n$. 

In \cite{BH} Baraglia and Huang show that the two Special K\"ahler structures -- the parameter space of fibres of $h:{\mathcal M} \rightarrow {\mathcal B}$ and of curves $S\subset T^*\Sigma$ -- are the same. We now give a slightly different account of this result.

\section{Two families of Lagrangians}
\subsection{Fibres of the integrable system} 
Consider again the tangent space to the moduli space of Higgs bundles ${\mathcal M}$. As noted above this is defined by $(\dot A,\dot \Phi)$ 
satisfying  $\bar\partial_A\dot\Phi+[\dot A,\Phi]=0$ modulo  $(\dot A,\dot\Phi)=(\bar\partial_A\psi,-[\psi,\Phi])$. 
This space has a holomorphic interpretation as  the  hypercohomology $\K^1$ of the complex of sheaves: 
$${\mathcal O}(\End V)\stackrel {[\Phi, -]}\rightarrow {\mathcal O}(\End V\otimes K).$$
One of the two spectral sequences associated   to this complex gives, using the vanishing theorem coming from stability, the following:
\begin{equation}
0 \rightarrow H^1(\Sigma, \ker \Phi)\rightarrow \K^1\rightarrow H^0(\Sigma, \coker \Phi)\rightarrow 0.
\label{tangent}
\end{equation}
Here $\ker\Phi$  is the sheaf of centralizers of $\Phi$ and $\coker \Phi$ the cokernel sheaf of the adjoint action of $\Phi$. Since $\bar\partial_A\dot\Phi=-[\dot A,\Phi]$ we see that $\dot\Phi$ is holomorphic modulo the image of $\ad \Phi$ and so defines the image of  $(\dot A,\dot\Phi)$ in  $H^0(\Sigma, \coker \Phi)$.

If the spectral curve is smooth then $\Phi$ is everywhere regular and hence the  centralizer is generated by $1,\Phi,\dots,\Phi^{n-1}$ where $\Phi^m:K^{-m}\rightarrow \End V$. It follows that 
  $\ker \Phi\cong {\mathcal O}\oplus K^{-1}\oplus \cdots\oplus K^{-(n-1)}$ and similarly 
$\coker \Phi\cong K\oplus K^{2}\oplus \cdots\oplus K^n.$

If $p: {\mathcal O}(\End V) \rightarrow \coker \Phi$ is the projection then $\tr (\Phi^m\dot\Phi)=\tr (\Phi^mp(\dot\Phi))$ and then the right hand arrow in (\ref{tangent}) is 
$$(\dot A,\dot\Phi)\mapsto (\tr\dot\Phi,\tr(\Phi\dot\Phi), \dots,\tr(\Phi^{n-1}\dot\Phi))\in \bigoplus_{m=1}^n H^0(\Sigma, K^m)$$
We recognize this as the derivative of the projection $h:{\mathcal M}\rightarrow {\mathcal B}$ of the integrable system, and so the tangent space of a fibre is  
$$H^1(\Sigma, \ker \Phi)\cong \bigoplus_{m=0}^{n-1} H^1(\Sigma, K^{-m}).$$
A tangent vector to the fibre is thus represented by 
\begin{equation}
\dot A=\sum_ {m=0}^{n-1} \alpha_m\Phi^m\in \Omega^{01}(\Sigma, \End V)
\label{dot}
\end{equation}
where 
$\alpha_m\in \Omega^{01}(\Sigma, K^{-m})$.

The symplectic form $\omega^c$ on ${\mathcal M}$ gives a nondegenerate pairing between the tangent bundle of the fibre and base. In terms of the representatives above this is 
\begin{equation}
\int_{\Sigma} \tr(\dot A\dot\Phi)=\sum_ {m=0}^{n-1}\int_{\Sigma} \alpha_m\tr(\Phi^m\dot\Phi)=\sum_ {m=0}^{n-1}\langle [\alpha_m], b_m\rangle
\label{pair}
\end{equation} 
where $(b_1,\dots, b_m) \in H^0(\Sigma, K)\oplus \cdots\oplus H^0(\Sigma, K^n)$ is a tangent vector to the base and $\langle a, b\rangle$ is Serre duality on $\Sigma$. 

The fibre is identified with the Jacobian of the spectral curve $S$ by the direct image.  Recall that $x$ is a section of $\pi^*K$ and $\Phi=\pi_*x$, so $\dot A$ in (\ref{dot}) defines naturally 
$$\sum_ {m=0}^{n-1} \pi^*[\alpha_m]x^m\in H^1(S,{\mathcal O})$$
a tangent vector to the Jacobian. As noted, the canonical bundle of $S$ is isomorphic to $\pi^*K^n$ but we need to examine the isomorphism more closely  to compare the Higgs bundle pairing above with the Serre duality pairing on the spectral curve $S$.

 \subsection{Curves in $T^*\Sigma$}\label{Tstar} 
 The other family of Lagrangians are the compact curves $S\subset \vert K\vert$ in the linear system $\pi^*K^n$, each given by the vanishing of a  section $s=x^n+a_1x^{n-1}+\cdots + a_n$. The surface $\vert K\vert=T^*\Sigma$ is the total space of the cotangent bundle of $\Sigma$, which has a canonical symplectic form $\omega$. If $z$ is a local coordinate on $\Sigma$ then $dz$ trivializes $K$ and $x\in \pi^*K$ is another local coordinate. The canonical 2-form is then $dx\wedge dz$.  
          
     The normal bundle $N$ of $S$ is $\pi^*K^n$ and the derivative of $s$ on $S$, assuming $S$ is smooth, is a nonvanishing    section of $    \pi^*K^n\otimes N^*$.     Since the canonical bundle of $T^*\Sigma$ is trivial, $N\cong K_S$ and this gives an isomorphism $K_S\cong \pi^*K^n$. The function $z$ is a local coordinate on $S$ if $s_x\ne 0$ and locally the derivative of $s$ is $s_xdx+s_zdz$. Then the isomorphism is locally 
     $$dz \mapsto s_x=nx^{n-1}+(n-1)a_1x^{n-2}+\cdots +a_{n-1}.$$
     
     The Abel-Jacobi map $S\rightarrow \Jac(S)$ realizes an isomorphism between holomorphic 1-forms on $\Sigma$ and $S$ and also between $H^1(\Sigma, \Z)$ and $H^1(\Jac(S),\Z)$. The derivative of this map takes a tangent vector at a point $p\in S$ to a tangent vector to $\Jac(S)$. It comes from the long exact sequence of 
     $$0\rightarrow {\mathcal O}_S\rightarrow {\mathcal O}_S(p)\rightarrow {\mathcal O}_p(p)\rightarrow 0.$$
 The derivative of the section of ${\mathcal O}(p)$ which vanishes at $p$ identifies the normal bundle of $p$ (i.e. the tangent space at $p$) with the line bundle ${\mathcal O}(p)$ so canonically $H^0(p, {\mathcal O}_p(p))\cong T_p$ and the connecting homomorphism  $\delta: T_p\cong H^0(p, {\mathcal O}_p(p))\rightarrow H^1(S,{\mathcal O})\cong T\!\Jac(S)$ gives the corresponding tangent vector.
  
If $\pi(p)$ is not a branch point then $x$ takes $n$ distinct values $\lambda_1,\dots,\lambda_n$ in a neighbourhood of $p$ so suppose $x=\lambda_1$ near $p$. Then define in a neighbourhood $$q(x)=\frac{(x-\lambda_2)\dots(x-\lambda_n)}{(\lambda_1-\lambda_2)\dots(\lambda_1-\lambda_n)}=c_0+c_1x+\cdots +c_{n-1}x^{n-1},$$
 equal to $1$ when $x=\lambda_1$ and zero for $\lambda_i, i>1$.  Note that the denominator is $s_x(\lambda_1)$, the identification of $\pi^*K^n$ with $K_S$. Applying $\delta$ gives a class in  $H^1(S,{\mathcal O})$ of the form 
 $$\sum_ {m=0}^{n-1} \pi^*\delta(c_m)x^m\in H^1(S,{\mathcal O}).$$
where $\delta (c_m)\in H^1(\Sigma, K^{-m})$. 

Serre duality on $S$ pairs this localized class with a holomorphic 1-form by evaluation at $p$. Let $b_i$ be a holomorphic section of $K^i$, and consider $\pi^*b_ix^{n-i}$ a section of $\pi^*K^n$. Then $\pi^*b_ix^{n-i}q(x)$ modulo $(x-\lambda_1)\dots(x-\lambda_n)$ is equal to $b_i(\pi(p))\lambda_1^{n-i}/s_x(\lambda_1)$ at $p$  since $x^{n-i}-\lambda^{n-i}$ is divisible by $x-\lambda_1$.  

On the other hand, the pairing via Higgs bundles in (\ref{pair}) is $c_ib_i(\pi(p))$ where $c_i=(-1)^{n-i}\tau_{n-i}/s_x(\lambda_1)$ and $\tau_k$ is the $k$-th elementary symmetric function in $\lambda_2,\lambda_3,\dots,\lambda_n$.

In terms of the full symmetric functions $\sigma_i$ we have 
 $\tau_k=\sigma_k-\lambda_1\sigma_{k-1}+\lambda_1^2\sigma_{k-2}+\cdots + (-1)^k\lambda_1^k$ so the homomorphism
\begin{equation}
(b_1,b_2,\dots, b_n) \mapsto \sum_{j=1}^n\pi^*b_j(x^{n-j}+\pi^*a_1x^{n-j-1}+\cdots + \pi^*a_{n-j})
\label{biso}
\end{equation}
identifies the tangent space of ${\mathcal B}$ with $H^0(S,\pi^*K^n)\cong H^0(S,K_S)$ in such a way that the Higgs bundle pairing is equal to the Serre duality pairing. 
This is Proposition 3.3 in \cite{B}.

\subsection{The Special K\"ahler metric}\label{SK}
Let ${\mathcal B}^{\mathrm{reg}}\subset {\mathcal B}$  denote the space of smooth spectral curves, then we have identified the family of Jacobians  with the Higgs bundle moduli space in a form where  the pairing given by the symplectic forms coincides.  The Abel-Jacobi map identifies holomorphic one-forms and their periods on the spectral curve with those of the Jacobian and it follows that the two Special K\"ahler metrics are the same.

Now the $\C^*$-action on ${\mathcal M}$ is given by $\Phi\mapsto \lambda\Phi$ and $\det (x-\lambda \Phi)=\lambda^n\det(\lambda^{-1}x-\Phi)$. So the action on spectral curves $S\subset T^*\Sigma$  is rescaling $x\mapsto \lambda^{-1}x$ of the fibre of the cotangent bundle. Call the associated vector field $Y$, then   $i_Y\omega=\theta$ is the canonical 1-form on the cotangent bundle, and we have the  following: 
\begin{prp} \label{potential}(\cite{BH}) A K\"ahler potential for the Special K\"ahler metric on ${\mathcal B}^{{\mathrm {reg}}}$ is given by 
$$K=\Imag \frac{1}{4}\int_S\theta\wedge\bar\theta.$$
\end{prp}

It is instructive to see the relationship between the two holomorphic 1-forms explicitly. In our local coordinates on $T^*\Sigma$ we have $\theta=xdz$. Using the identification of $K_S$ with $\pi^*K^n$ this becomes 
$$xs_x=x(nx^{n-1}+(n-1)a_1x^{n-2}+\cdots +a_{n-1})=n(-(a_1x^{n-1}+\cdots +a_n))+(n-1)a_1x^{n-1}+\cdots $$
$$= -(a_1x^{n-1}+2a_2x^{n-2}+\cdots  + na_n)$$
using the equation of the spectral curve. 

The $\C^*$-action on the Higgs field gives $\dot\Phi=\Phi$ and so $\tr(\Phi^m\dot\Phi)=\tr\Phi^{m+1}$ and  we substitute $b_j=\tr(\Phi^j)$ in equation (\ref{biso}). Newton's identity for sums of powers and elementary symmetric functions gives
$$\sum_{j=1}^kb_ja_{k-j}=-ka_k$$
and using this in (\ref{biso}) gives the above formula.  
\section{Singular curves}
\subsection{Singular curves and singular fibres}
If $S\subset T^*\Sigma$ is now a singular curve then the formula of Proposition \ref{potential} still makes sense and a natural question is to understand the significance of the corresponding potential. When $S$ is singular, the fibre of the integrable system becomes singular. This means that there are points where some of the Hamiltonian vector fields corresponding to functions on the base vanish, or equivalently they are critical points of those functions. More geometrically, it is where  the derivative of $h:{\mathcal M}\rightarrow {\mathcal B}$ is not surjective. 

We observed that if $\Phi$ is regular everywhere then surjectivity holds so 
suppose $\Phi$  is regular except on a divisor $D\subset \Sigma$ (for example if the spectral curve  $S$ is reduced), then the  powers of $\Phi$ still lie in the centralizer but now there is a quotient sheaf ${\mathcal C}$ supported on $D$ 
$$0\rightarrow {\mathcal O}\oplus K^{-1}\oplus \cdots\oplus K^{-(n-1)}\rightarrow \ker\Phi\rightarrow {\mathcal C}\rightarrow 0.$$

From the long exact sequence 
$$\rightarrow H^0(D,{\mathcal C})\stackrel{\delta}\rightarrow H^1(\Sigma, {\mathcal O}\oplus K^{-1}\oplus \cdots\oplus K^{-(n-1)})\rightarrow $$
sections of ${\mathcal C}$ define by Serre duality  linear functions on the base of the integrable system ${\mathcal B} = H^0(\Sigma,K)\oplus \cdots\oplus H^0(\Sigma, K^n).$
These are functions which are critical at the point $(V,\Phi)\in {\mathcal M}$. 

In rank $2$, this critical locus was discussed in \cite{Hit2}, with emphasis on a notion of nondegeneracy of the critical locus. Here we make the simplifying assumption that the singularities of $S$ are {\it ordinary singularities} in the classical terminology of Walker \cite{W}, namely that each singularity of multiplicity $m$ has $m$ distinct tangents. When $m=2$ this is a node. If $S$ has $k$ ordinary singularities of multiplicity $m_j$ and $\tilde S$ is the normalization of $S$ then its genus is 
\begin{equation}
 g(\tilde S)= g(S)-\sum_1^k\frac{m_j(m_j-1)}{2}.
 \label{genus}
 \end{equation} 
The fibres of the Higgs bundle integrable system are described in general by torsion-free sheaves on the spectral curve which can be regarded as the direct image of line bundles on a partial normalization. So consider the case of an ordinary singularity of $S$. At this point $m$ eigenvalues of $\Phi$ coincide. For simplicity assume that these are zero  and $z$ is a local coordinate on $\Sigma$, then the equation of $S$ is locally of the form $(x-\lambda_1z)(x-\lambda_2 z)\dots(x-\lambda_mz) +$ higher order terms. Distinct tangents mean that the $\lambda_i$ are distinct. 

The normalization of a neighbourhood of the singularity consists of $m$ disjoint discs $U_i$ on which $x=\lambda_i z+\cdots$ and so the direct image of $x$ contributes a component of the Higgs field of the form $z\phi + \cdots$ where $\phi$ is the diagonal matrix with entries $\lambda_1,\dots,\lambda_m$. The quotient sheaf ${\mathcal C}$ above is thus generated by $\phi,\phi^2,\dots,\phi^{m-1}$ modulo powers of $z\phi$ or 
$\sum_{j=1}^{m-1} p_k(z)\phi^k$
where $p_k(z)$ is a polynomial of degree $k-1$. Then at this singularity $\dim H^0(D,{\mathcal C})= 1+2 +\cdots +(m-1)=m(m-1)/2$. 

Let ${\mathcal B}_d\subset {\mathcal B}$ denote the subvariety of spectral curves where $d$ denotes the singularity type ($k$ ordinary singularities together with multiplicities $m_k$). Note that ${\mathcal B}_d$ is preserved by the $\C^*$-action. A tangent vector at a smooth point is annihilated by   $\delta H^0(D,{\mathcal C})$ whose dimension is $\sum m_k(m_k-1)/2$ and so the tangent space of ${\mathcal B}_d$ has dimension $g(\tilde S)$ from (\ref{genus}).  The generic fibre is isomorphic to the Jacobian of $\tilde S$ and so we have a subvariety of dimension $2g(\tilde S)$. 

This is a symplectic submanifold and a subintegrable system ${\mathcal M}_d$, in fact the local structure can be realized as a symplectic quotient. If  ${\mathcal B}_d$ is locally defined by functions $(f_1,f_2,\dots, f_m)=0$, then these constitute the moment map for  a $\C^m$-action which translates the torus fibres $\cong \Jac(\tilde S)$. The symplectic quotient fibres over ${\mathcal B}_d$ and is locally holomorphically equivalent to the symplectic submanifold. As an integrable system in its own right, ${\mathcal B}_d$ inherits a Special K\"ahler structure and we shall determine next the K\"ahler potential by using the canonical 1-form $i_X\omega^c$.

\subsection{The subintegrable system} 
The symplectic pairing on ${\mathcal M}_d$  provides an isomorphism between the tangent space of a fibre, namely $H^1(\tilde S,{\mathcal O})$, and the cotangent space of the fibre ${\mathcal B}_d\subset {\mathcal B}$. In particular there is an  identification of  the sections of $\pi^*K^n$ on $S$ which are annihilated by $\delta H^0(D,{\mathcal C})$ with holomorphic one-forms on $\tilde S$, the normalization of $S$. We need to see this concretely.

Near an $m$-fold singularity $x=0=z$, the section of $\pi^*K^n$ defining $S$ is of the form 
$$s=(x-\lambda_1z)(x-\lambda_2z)\dots(x-\lambda_mz) + \cdots$$
 Then on the neighbourhood $U_i\subset \tilde S$ on which $x=\lambda_1z+\cdots$ we have 
$$s_x=z^{m-1}dx(\lambda_1-\lambda_2)\dots(\lambda_1-\lambda_m)+\cdots$$ 
and since the $\lambda_i$ are distinct this gives an isomorphism to $K_{\tilde S}$ after dividing by $z^{m-1}$. So if $p_i\in U_i$ are the points in the resolution of the singularity the derivative of $s$ gives an isomorphism  $K_{\tilde S}\cong \pi^*K^n(-(m-1)(p_1+\cdots + p_m))$. (This reduces the degree by $m(m-1)$ providing another argument for the genus $g(\tilde S)$ whose degree is $2g(\tilde S)-2$.) 

Now consider $(b_1,b_2,\dots,b_{n}) \in H^0(\Sigma,  K\oplus K^2\oplus \cdots \oplus K^n)$ annihilated by $\delta H^0(D,{\mathcal C})$. A class $\delta(p_k)\in H^1(\Sigma, K^{-k})$ is represented by  $c_0+c_1z +\dots +c_kz^{k-1}$ modulo $z^k$ and if all such classes annihilate a section $b_k$ of $K^k$ then $b_k$ is divisible by $z^k$.  Moreover since the equation of $S$ is $x^n+a_1x^{n-1}+\cdots= (x-\lambda_1z)(x-\lambda_2 z)\cdots(x-\lambda_mz) + \cdots$, $a_i$ is divisible by $z^i$ if $i\le m$.  Finally, on the open set  $U_i$, $x= \lambda_i z +$ higher order terms. Putting these three facts together,  Equation \ref{biso} shows that the isomorphism from $\pi^*K^n$ to the canonical bundle given by the derivative of the section $s$ at smooth points of $S$ extends to the normalization $\tilde S$. 

We can now follow the same procedure as in Sections \ref{Tstar} and \ref{SK} by choosing a generic point $p\in S$  to identify the periods of the form $i_X\omega^c$ on ${\mathcal M}$ restricted to ${\mathcal M}_d$ and the canonical one-form on $T^*\Sigma$ pulled back under the normalization map $p:\tilde S\rightarrow S\subset T^*\Sigma$. 

To complete the picture of the Special K\"ahler metric we  observe that the K\"ahler class on ${\mathcal M}$ was the first Chern class of the determinant line bundle. If we normalize the degree of $V$ to be $n(g-1)$ then there is a determinant divisor -- the subvariety   on which $H^0(\Sigma, V)\ne 0$. But if $V$ is the direct image of a line bundle  $ L$ on $\tilde S$ then $H^0(\Sigma, V)=H^0(\tilde S, L)$ by definition, so this class on the Jacobian is the class of the theta-divisor. The Special K\"ahler structure on ${\mathcal B}_d$ is therefore determined via the geometry of the normalization of the spectral curve. In particular, since $S$ and $\tilde S$ differ on a finite set of points, we  obtain: 

\begin{prp} \label{potentialsing}  A K\"ahler potential for the Special K\"ahler metric on ${\mathcal B}_d$ is given by 
$$K=\Imag \frac{1}{4}\int_S\theta\wedge\bar\theta.$$
\end{prp}
\section{Features and examples}
\subsection{Translation invariance}
Let $\alpha$ be a holomorphic 1-form on $\Sigma$ then it defines a diffeomorphism of $T^*\Sigma$ by translation in the fibre  direction, mapping the spectral curve $x^n+a_1x^{n-1}+\cdots + a_n$ to $(x-\alpha)^n+a_1(x-\alpha)^{n-1}+\cdots+a_n$. The K\"ahler potential $K$ is now 
$$K=\Imag \frac{1}{4}\int_S\left(\theta\wedge\bar\theta+\pi^*\alpha\wedge\bar\theta+\theta\wedge\pi^*\bar\alpha+\pi^*\alpha\wedge\pi^*\bar\alpha\right).$$
This is the addition of $f+\bar f$ where $f$ is holomorphic on ${\mathcal B}$ and so leaves the metric unchanged. Thus $H^0(\Sigma, K)$ acts as an isometric group of translations. The quotient can be identified with the Special K\"ahler manifold for the moduli space of $PGL(n,\C)$ Higgs bundles. 

\subsection{Rank $2$ bundles}
The lowest dimensional case for providing an example is when $V$ has rank $2$ with $\Lambda^2V$ trivial and $\tr\Phi=0$ -- these are $SL(2,\C)$-Higgs bundles and the fibre of the integrable system is now the Prym variety inside the Jacobian of $S$. The spectral curve  has a simple formula $x^2-q=0$ where $q\in H^0(\Sigma, K^2)$. In this case, since $x=\pm \sqrt{q}$,  the formula for the K\"ahler potential can be written as an integral over $\Sigma$
$$K=\frac{1}{2}\int_{\Sigma} \sqrt{q\bar q}.$$
Here $q\bar q$ is a section of $K^2\bar K^2$ so that $\sqrt{q\bar q}$ is a positive section of $K\bar K$ and therefore a measure. 
This is discussed also in \cite{MSWW}. The singularity type to determine ${\mathcal B}_d$ is now just a question of the multiple zeros of the quadratic differential $q$ and ordinary singularities are nodes.

Even simpler is the case where $g=2$ and $\Sigma$ is the hyperelliptic curve 
$y^2=f(z)=(z-z_1)\dots(z-z_6)$
for then a quadratic differential can be written as 
$$q=(c_0+c_1z+c_2z^2)\frac{dz^2}{f(z)}$$
where $(c_0,c_1,c_2)\in \C^3\cong {\mathcal B}.$
The K\"ahler potential is then an integral over the complex plane $\C$:
$$K=\frac{1}{2}\int_{\C} \frac{\vert c_0+c_1z+c_2z^2\vert}{\vert f(z)\vert}dzd\bar z.$$

For $q$ to have one double zero requires $c_0+c_1z+c_2z^2$ to vanish at a branch point $z=z_i$. There are six of these and so six components of ${\mathcal B}_1$. The quadratic differential is  $q=(a_0+a_1z)(z-z_i)dz^2/f(z)$  where $a_0+a_1z_j\ne 0$ ensures that there are no more singularities. Each component is thus the complement of  six lines through the origin of $\C^2$.

Two double zeros implies $q=a(z-z_i)(z-z_j)dz^2/f(z)$ for $i\ne j$ so ${\mathcal B}_2$ has $15$ components each one a copy of $\C^*$. Here the K\"ahler potential is $K=c\vert a\vert$ where $c$ is the constant
$$c=\frac{1}{2}\int_{\C} \frac{\vert (z-z_i)(z-z_j)\vert}{\vert f(z)\vert}dzd\bar z.$$
Setting $a=w^2$ this gives $K=cw\bar w$ which is a flat metric on $\C^*$. 
The Special K\"ahler structure is non-trivial however since $\nabla(dw)=0$ implies that the connection $\nabla$ in the coordinate $a$ is 
$$\nabla=\frac{d}{da}+\frac{1}{2a}$$ and has a pole at $a=0$ with nontrivial holonomy. 

\begin{rmk}
In general we observe that the complex structure on ${\mathcal B}^{\mathrm{reg}}$ extends to ${\mathcal B}_d\subset {\mathcal B}$ and Proposition \ref{potentialsing} shows that the K\"ahler potential extends,  continously but  not smoothly. It is the flat connection component of the Special K\"ahler structure which acquires a singularity, but this is typically a logarithmic singularity which induces a connection on the divisor.  For the Gauss-Manin connection considered here the induced  connection is presumably the flat connection for the Special K\"ahler structure on the divisor ${\mathcal B}_1$ where the spectral curve acquires a node. 
\end{rmk}
\subsection{Fixed points of a $\Z_n$-action} \label{fix}
The last example admits a generalization to further cases where the metric is flat. Consider the action on ${\mathcal M}$ of tensoring $V$ with a line bundle $U$ such that $U^n$ is trivial, where $n=\rk V$ as usual.  The fixed point set of this $\Z_n$-action is described in \cite{FGOP}.

Firstly, having a fixed point means we have an isomorphism $\varphi:V\rightarrow V\otimes U$. This is a Higgs bundle where the twist is by $U$ instead of $K$ but the general theory tells us that $V$ is the direct image of a line bundle $L$ on the unramified cyclic $n$-fold cover $p:C\rightarrow \Sigma$ associated to $[U]\in H^1(\Sigma, \Z_n)$. Let $\sigma$ denote a generator of the $\Z_n$-action on $C$. 

Since $p^*U$ is trivial, $p^*\varphi$ is an endomorphism of $p^*V$. To be a fixed point the Higgs field $\Phi$ must commute with $\varphi$  and so on $C$ 
$$p^*\Phi = c_0+c_1p^*\varphi+\cdots + c_{n-1}p^*\varphi^{n-1}$$
where $c_i$ are sections of $p^*K$. Then 
$\Phi = c_0+c_1\varphi+\cdots +c_{n-1}\varphi^{n-1}$
for sections $c_k\in H^0(\Sigma, KU^{-k})$. 

Since $p:C\rightarrow \Sigma$ is unramified, $p^*K\cong K_C$ and then the $c_i$ are  sections of $K_C$ which transform as $\omega^{-i}$ under the action of $\sigma$, and an arbitrary 1-form on $C$ can be decomposed this way.  The data is then just a holomorphic 1-form and a line bundle on $C$. 
This gives an integrable subsystem  with a flat Special K\"ahler metric.

Since  $(p^*\varphi)^n=1$ its eigenvalues are $1,\omega,\omega^2, \dots$ where $\omega=e^{2\pi i/n}$ and so the  eigenvalues of $p^*\Phi$ are $c_0+c_1+\dots +c_{n-1}, c_0+\omega c_1+\dots +\omega^{n-1}c_{n-1}$, etc. Put another way, $\gamma = c_0+c_1+\dots +c_{n-1}$ is a 1-form on $C$ and the eigenvalues are $\gamma,\sigma^*\gamma, (\sigma^2)^*\gamma,\dots$. Hence the spectral curve $S$ has equation 
\begin{equation}
\prod_{k=0}^{n-1} (x-(\sigma^k)^*\gamma)=0.
\label{product}
\end{equation}
On $C$ we have $n$ single-valued roots of the characteristic polynomial so there is a map from $C$ to the spectral curve $S$. Both cover $\Sigma$ by taking eigenvalues of $\Phi$ so the map is generically one-to-one. Since $C$ is nonsingular it is the normalization of $S$.

If $\gamma$ has a simple zero at $c\in C$ then the other eigenvalues are 1-forms with zeros at $ \sigma(c), \sigma^2(c),\dots$. Then equation (\ref{product}) shows that $S$ has an ordinary singularity of multiplicity $n$   over $p(c)\in \Sigma$. 
Since $C$ is an unramified $n$-fold covering of $\Sigma$, we have  $2g(C)-2=n(2g-2)$ and $\gamma$ has $n(2g-2)$ zeros. Thus the number of singularities of this type on $S$ is $(2g-2)$. But the arithmetic genus of $S$ is $n^2(g-1)+1$ and 
$$n^2(g-1)+1- (2g-2)\frac{n(n-1)}{2}=n(g-1)+1=g(C)$$
so $S$ has no more singularities. 

It follows that the singularity type for this subintegrable system consists of $n$ ordinary $n$-fold singularities.

The simplest case is the group $SL(2,\C)$  where $U\in H^1(\Sigma, \Z_2)$ is a non-trivial line bundle of order $2$. Then $H^0(\Sigma, KU)$ has dimension $g-1$  and a section $s$ defines the quadratic differential $q=s^2$. Then the formula (\ref{potentialsing}) gives a K\"ahler potential 
$$K=\frac{1}{2}\int_{\Sigma} \sqrt{q\bar q}=\frac{1}{2}\int_{\Sigma}s\bar s.$$
This Hermitian expression in $s$ gives the flat metric on $H^0(\Sigma, KU)\backslash \{0\}/\pm 1$. 
\subsection{Hyperk\"ahler metrics} 
Any Special K\"ahler manifold gives a semiflat hyperk\"ahler metric on a manifold modelled on its cotangent bundle \cite{DF},\cite{Hit1}. The moduli space of Higgs bundles ${\mathcal M}$ has a natural $L^2$-metric which has been considered in more detail quite recently as in \cite{MSWW}, \cite{LF}. In particular, over ${\mathcal B}^{\mathrm{reg}}$ it is shown that the semiflat metric is a good approximation to the actual metric. It seems likely, then,  that  the formula here for the potential could provide an approximation on the lower dimensional loci ${\mathcal M}_d$. 

In fact, in support of this suggestion, the fixed point set in  Section \ref{fix} is a hyperk\"ahler submanifold of the genuine metric on ${\mathcal M}$ and (because one can take the  direct image of a flat connection in an unramified cover) the induced metric  is also flat. 

\subsection*{Acknowledgements}
The author wishes to thank  the EPSRC Programme Grant ``Symmetries and Correspondences"    and the Instituto de Ciencias Matem\'aticas, Madrid  for support during the preparation of this work, and David Baraglia for useful discussions.  This paper is based on a talk given at the conference ``Between topology and quantum field theory" at the University of Texas in Austin in celebration of Dan Freed's 60th birthday.

\end{document}